\documentclass{article}[11pt]
\usepackage{geometry}
\usepackage{amsmath, amssymb}
\usepackage{hyperref,cleveref}
\usepackage{subcaption}
\Crefname{equation}{Eq.}{Eqs.}
\Crefname{figure}{Fig.}{Figs.}
\Crefname{table}{Tab.}{Tabs.}

\usepackage{changes}

\usepackage{amsthm}
\usepackage{bm}
\usepackage{xcolor}
\usepackage{graphicx}
\usepackage{authblk}
\usepackage{mathtools}
\usepackage{todonotes}
\usepackage{tikz}
\usetikzlibrary{positioning,calc}
\usepackage{booktabs}
\newcommand{\rT}[1]{#1^{\textsf{T}}}
\newcommand{\fA}{\mathbf{A}}
\newcommand{\fa}{\mathbf{a}}
\newcommand{\ff}{\mathbf{f}}
\newcommand{\fx}{\mathbf{x}}

\newcommand{\fz}{\mathbf{z}}
\newcommand{\fU}{\mathbf{U}}
\newcommand{\fV}{\mathbf{V}}
\newcommand{\fZ}{\mathbf{Z}}
\newcommand{\fSigma}{\mathbf{\Sigma}}

\usepackage{appendix}
 
\appendixtitleon
\appendixtitletocon

\newcommand{\R}{\mathbb{R}}
\newcommand{\N}{\mathbb{N}}

\newcommand{\damkoehlerNumber}{\mathrm{Da}_L}
\newcommand{\pecletNumber}{\mathrm{Pe}_L}

\newcommand{\train}[1]{#1_{\mathrm{train}}}
\newcommand{\test}[1]{#1_{\mathrm{test}}}
\newcommand{\erel}{e_{\mathrm{rel}}}
\newcommand{\featureMap}{\Phi}
\newcommand{\featureMapPCA}{\featureMap_{\mathrm{PCA}}}

\newcommand{\numFeatures}{n_{\mathrm{f}}}
\newcommand{\numTimesteps}{n_t}
\newcommand{\numSamples}{n_{\mathrm{s}}}
\newcommand{\numSamplesTrain}{n_{\mathrm{s,train}}}
\newcommand{\numSamplesTest}{n_{\mathrm{s,test}}}
\newcommand{\numGeom}{2 N_V}

\begin{document}
\title{Greedy Kernel Methods for Approximating Breakthrough Curves for Reactive Flow from 3D Porous Geometry Data}
\author[1]{Robin Herkert \thanks{robin.herkert@mathematik.uni-stuttgart.de}}
\author[1,3]{Patrick Buchfink
	\thanks{p.buchfink@utwente.nl}}
\author[1,4]{Tizian Wenzel \thanks{tizian.wenzel@uni-hamburg.de}}
\author[1]{Bernard Haasdonk \thanks{haasdonk@mathematik.uni-stuttgart.de}}
\author[2]{Pavel Toktaliev \thanks{pavel.toktaliev@itwm.fraunhofer.de}}
\author[2]{Oleg Iliev \thanks{oleg.iliev@itwm.fraunhofer.de}}

\affil[1]{Institute for Applied Analysis and Numerical Simulation, University of Stuttgart, Stuttgart, Germany}
\affil[2]{Fraunhofer ITWM, Kaiserslautern, Germany, Technical University Kaiserslautern, Kaiserslautern, Germany}
\affil[3]{Department of Applied Mathematics, University of Twente, Enschede, The Netherlands}
\affil[4]{Department of Mathematics, Universität Hamburg, Hamburg, Germany}

\maketitle %

\abstract{We address the challenging application of 3D pore scale reactive flow under varying geometry parameters. The task is to predict time-dependent integral quantities, i.e.,\ breakthrough curves, from the given geometries. As the 3D reactive flow simulation is highly complex and computationally expensive, we are interested in data-based surrogates that can give a rapid prediction of the target quantities of interest. This setting is an example of an application with scarce data, i.e.,\ only having available few data samples, while the input and output dimensions are high. In this scarce data setting, standard machine learning methods are likely to fail. Therefore, we resort to greedy kernel approximation schemes that have shown to be efficient meshless approximation techniques for multivariate functions. We demonstrate that such methods can efficiently be used in the high-dimensional input/output case under scarce data. Especially, we show that the vectorial kernel orthogonal greedy approximation (VKOGA) procedure with a data-adapted two-layer kernel yields excellent predictors for learning from 3D geometry voxel data via both morphological descriptors or principal component analysis.}\\
\\
\noindent
\textbf{Keywords:}
Machine learning, Kernel methods, Two-layered Kernels, Porous media, Breakthrough curves

\section{Introduction}
Reactive flow in porous media plays an important role for many industrial, environmental and biomedical applications. Since the reactions occur at the pore scale 3D pore scale simulations are very important. At the same time, pore scale measurements are difficult or impossible, and usually some averaged quantities are measured. 
Such quantity of interest, which can be measured, is the breakthrough curve, i.e.,\ the time -dependent integral of the species over the outlet. These breakthrough curves can be computed from a reaction-advection-diffusion equation on a porous medium, which is numerically solved. However, this usually leads to very high computational costs which might be prohibitively high in a multiquery application, e.g., optimization of the geometry of the porous medium. In that case a surrogate model for the full order simulation model (FOM) is required. 
A promising way of obtaining an adequate surrogate model is the use of machine learning (ML) techniques. ML algorithms are developed and tested for porous media flows with different complexity. The literature in this area is very rich, and we exemplary cite some papers, just to point to the place of our research. Numerous papers discuss predicting permeability from microscale geometry. For example see \cite{Wu2018, marcato2022from, gaerttner2023estimating} 
and references therein.
At the next level of complexity, e.g., in \cite{Santos2020} the capability of using deep learning techniques has been shown in the case where the velocity field is predicted from the morphology of a porous medium. On the same topic, improvement is achieved by incorporating coarse velocities in the learning process, see \cite{Zhou22}.
In the last decade the application of ML techniques for simulation of passive and reactive transport rapidly grows. One of the directions in this case is using ML to predict reaction rates when those are very expensive in the case of complex reactions. This can be done without taking into account the geometry, see e.g. \cite{doeppel2024goal, Laloy2022speeding, silva2022rapid}, or predicting the rate from  structure features, see, e.g., \cite{Liu2022, marcato2022from, marcato2023phd}, to name just a few.
Recently, the potential of using machine learning models as surrogates for predicting  breakthrough curves for varying physical parameters, i.e.,\ Damk\"ohler and Peclet numbers, on a fixed porous medium geometry \cite{fokina2022performance,fokina2023machine} has been reported. 

In the current paper, we address the task of predicting the breakthrough curves for varying geometries of the porous medium with fixed Damk\"ohler and Peclet numbers. Our work is in the same direction as that of \cite{Liu2022, marcato2022from}, but there are essential differences. In \cite{Liu2022} ML addresses the impact of the structural features on the effective reaction rate to overcome the limitation of the well-mixing assumption. Pore scale reactive flow in an inert skeleton is considered there. In \cite{marcato2022from} pore scale colloid transport is considered as a part of a filtration problem. Steady-state problems are solved. In our case our pore scale geometry is in fact a two scale media, as the active washcoat particles are nanoporous. We consider the diffusion of the species within the washcoat, where the reaction occurs. Furthermore, while the other papers discuss different neural network algorithms, we consider a kernel-based method, namely, VKOGA. To the best of our knowledge, such studies have not yet been discussed in the literature.

Because of the computationally demanding simulation of the FOM, we are limited to a scarce data regime as we can only afford to compute a few samples. 
Furthermore, the input dimension (the number of elements of the  discretization of the porous medium) and the output dimension (the number of time steps used during the FOM simulation) are very high. 
This yields a challenging task for machine learning techniques \cite{Grigo19}. 
For this purpose, greedy kernel approximation schemes \cite{wirtz2013vectorial, Santin2021VKOGA, wenzel2023analysis} have been shown to be efficient.
These meshless approximation techniques can be combined with deep learning techniques to yield two-layered kernels \cite{wenzel2024data, wenzel2024application},
which have also already been successfully applied for certified and adaptive surrogate modeling for heat conduction \cite{wenzel2024application} as well as surrogate modeling of chemical kinetics  \cite{doeppel2024goal}.

Our work is organized as follows: In \Cref{sect:problem}, we introduce the underlying equations of the 3D porous medium reaction-advection-diffusion model.
In \Cref{sect:kernels}, we give an overview on greedy kernel methods and two-layered kernels. 
Numerical experiments are provided in \Cref{sect:num}
and we conclude our work in \Cref{sect:concl}.

\section{Problem Formulation}\label{sect:problem}
\subsection{Geometry}\label{subsect:geometry}
In the current paper we consider artificially generated voxel-based geometries on the unit cube $\Omega = \left[0,1\right]^3$ as geometries for porous filter fragments.
These consist of the solid skeleton, $\Omega_s$, free pores, $\Omega_f$, and effective porous space (washcoat or unresolved porosity region) $\Omega_w$. 
Accordingly, we assume that the porous geometry sample can be represented as a union of non-overlapping domains $\Omega = \Omega_s \cup \Omega_f \cup \Omega_w$. Each of the domains $\Omega_w, \Omega_f$ is represented as a union of non-overlapping volume voxels:
\begin{align*}
V_{i,j,k} = [(i - 1)h; ih) \times [(j - 1)h; jh) \times [(k - 1)h; kh), i,j,k = 1,..., N_h, \; h = 1/N_h,    
\end{align*}
where $N_h$ is the number of voxels in each dimension. 
This leads to a uniform grid with $N_V = 150^3$ elements (voxels). 
The assignment of the voxels to the three subdomains is stored within two boolean arrays
(a third one is not necessary as it can be calculated from the other two). 
That means each porous sample is represented by a vector $\fz \in \{0,1\}^{2 N_V}$.

\begin{figure}[t]
    \centering
    \begin{minipage}{0.5\textwidth}
        \includegraphics[width=\textwidth]{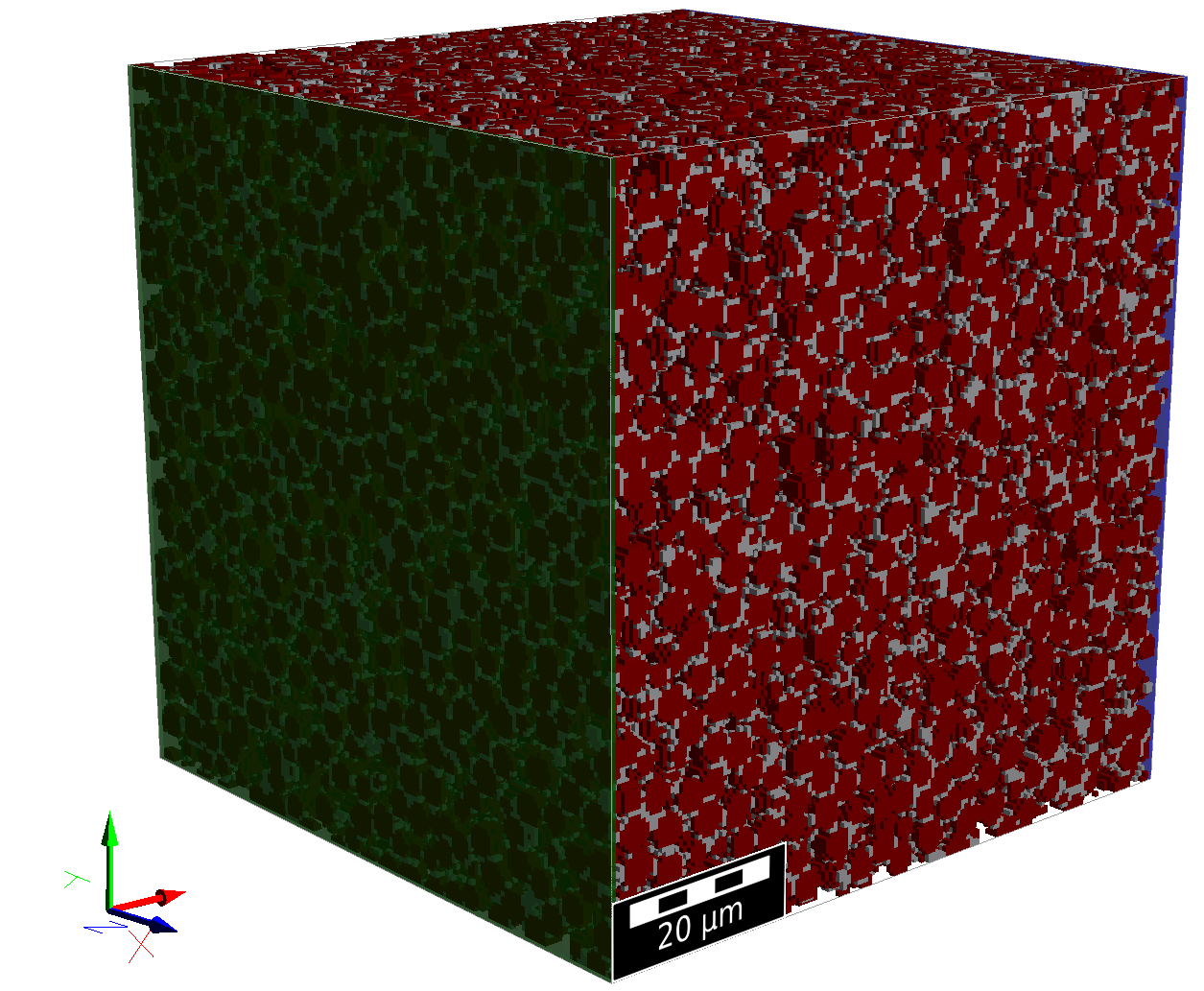}
    \end{minipage}
    \hspace{1cm}
    \begin{minipage}{0.35\textwidth}
        \includegraphics[width=\textwidth]{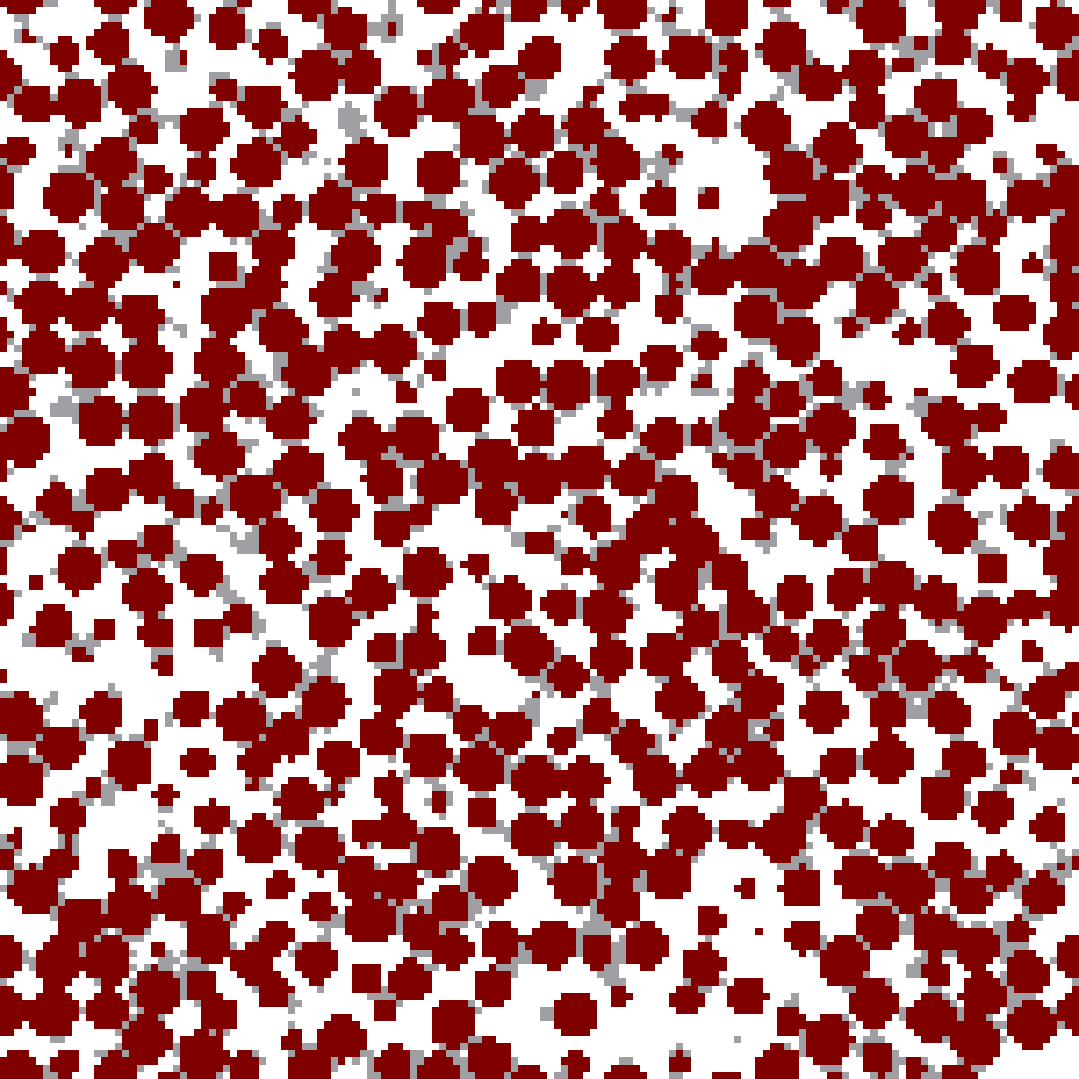}
    \end{minipage}
    \caption{Isometric (left) and middle plane (right) view of a typical porous sample; colors: brown - washcoat ($\Omega_w$), grey - solid (binder, $\Omega_s$), white (transparent) - free pores ($\Omega_f$), green - inlet boundary section, opposite of the (non-visible) - outlet section. GeoDict visualization \cite{Geodict2022}.}
    \label{fig:porousSample}
\end{figure}

Depending on the application and the production process, the structure of the porous space in real filters could be strongly irregular and highly anisotropic. 
Nevertheless, for certain types of filters, the micro structure of porous media could be represented as a combination of regular shaped nanoporous granules and solid binder material \cite{Kato2017}. 
We use this representation to reconstruct the porous domain $\Omega_f$ and the washcoat domain $\Omega_w$ and later, using existing experimental data on porosity, 
to generate a series of artificial porous geometries, which could be used as models for real filter fragments. 
The size of each sample in voxels, $N_h=150$, is chosen according to a representative volume element (RVE) study on the one hand and the amount of computational work to generate enough data for training on the other hand.
In order to generate each porous sample we use an analytical spheres piling algorithm for the washcoat region, then voxelize analytical spheres,
and distribute voxelized solid material (binder) uniformly and randomly, covering the washcoat skeleton surface. 
The main varying parameter during the generation of the samples was the washcoat volume fraction $\epsilon_w$, 
thus, in order to generate each sample the target value of $\epsilon_w$ was set, but depending on geometry realization, the resulting value of $\epsilon_w$ for a generated sample could differ by more than 5\%. 
A typical porous sample generated with this two-step procedure with porosity $\epsilon = 0.553$ is depicted in \Cref{fig:porousSample}.

Besides typical porous geometry characterization parameters, such as porosity $\epsilon$, 
phase volume fraction (for washcoat) $\epsilon_w$, 
and phase specific surface area (for free pores) $A_S$, 
we also compute integral Minkowski parameters as morphological features. 
A complete list of parameters with their definitions for each porous sample can be found in the supplementary \Cref{sect:appendix_S}.
In total, $59$ artificial porous samples were generated for further processing. 
The geometry generation time per sample varied greatly depending on the porosity, 
from the minimum of $6$ wall-clock seconds (wcs) to the maximum of $222$ wcs, with the average of $111$ wcs for a system with 40 cores Intel Xeon CPU E5-2600 v2, 2.8 GHz. 
On the contrary, the time for evaluation of geometric parameters was almost the same for all samples and did not exceed $2$ wcs. 

Both generation and parameter evaluation phases for each geometry were performed with the GeoDict software \cite{Geodict2022} and Python scripts.

\subsection{Governing Equations for FOM}
\label{subsec:governing_eq_for_fom}

Different approaches can be used to describe the flow and transport of chemical species through the porous medium at the pore-scale.
Many of them are based on solving convection-diffusion-reaction (CDR) equations for the species transport and (Navier-)Stokes equations for the bulk flow. 
Thus, we consider as FOM solution in the current paper, the solution of the CDR equation for a species of interest and additionally assume, 
based on typical filtration application conditions, that (i) flow Reynolds numbers are sufficiently small $Re \ll 1$ and Stokes equations are valid for any sub-region of the pore space $\Omega_f$, 
(ii) the concentration of the species of interest, $c$, is much smaller than the bulk mixture concentration.
This justifies a one-way coupling approach when the solution of the flow equations can be decoupled from the species transport equation. 
We also restrict ourselves to scalar CDR equations with a linear source term. This type of source term can be used to describe a first order sorption-desorption process in porous media (Henry isotherm). 
According to the aforementioned assumptions, the CDR equation for species concentration, $c \ge 0$, can be written in dimensionless form as:
\begin{equation}\label{eq:CDR}
\partial_t c  - \Delta c +  \pecletNumber \; \nabla \cdot  \left( {\bf v} c \right)  + \damkoehlerNumber \; c = 0, \quad x\in \left( \Omega_f \cup \Omega_w \right), ~ t>0, \\
\end{equation}
where $\damkoehlerNumber,\pecletNumber$ are the parameters: $\damkoehlerNumber = \frac{k_R L}{D}$ is the piecewise-constant Damk\"ohler number, $\pecletNumber = \frac{u_{in}L}{D}$ is the Peclet number, $L, D > 0$ are the characteristic length and diffusion coefficient respectively, $u_{in}$ is the inlet velocity, $k_R(x) \ge 0$ is the reaction rate constant and reaction occurs only in the washcoat region:
\begin{equation*}\label{eq:reactionRateK}
    k_R(x) =
    \begin{cases}
        0, \quad \mathbf{x} \in \Omega_f, \\
        k_r > 0, \quad \mathbf{x} \in \Omega_w.
    \end{cases}
\end{equation*}
In \Cref{eq:CDR} the characteristic diffusion time was used as a time scale. 
Moreover, we fix the parameter values, $(\pecletNumber,\damkoehlerNumber)=(5,0.1)$, and, thus, consider only a convection-dominated regime for each geometry. 
For the numerical solution \Cref{eq:CDR} was complemented with zero initial conditions, Dirichlet conditions for the inlet boundary section and zero Neumann conditions for the outlet boundary section, 
and all other boundaries, see \Cref{fig:porousSample}.
Based on the solution of \Cref{eq:CDR} a quantity of interest, the breakthrough curve, can be evaluated as an integral concentration over the outlet section of the geometry, $ \Gamma_{\mathrm{outlet}}$:
\begin{equation}\label{eq:bcurve}
a(t) = \int\limits_{\Gamma_{\mathrm{outlet}}} c(x,t) ~ \mathrm{d}\sigma.
\end{equation}
For the velocity ${\bf v}$ within \Cref{eq:CDR}, the stationary Stokes equations in $\Omega_f$ were considered:
\begin{equation}\label{eq:Stokes}
\begin{split}
\mu \Delta {\bf v} = \nabla p, \\
\nabla \cdot {\bf v} = 0,
\end{split}
\end{equation}
where $p: \Omega_f \rightarrow \R$ is the pressure and $\mu \in \mathbb{R}_+$ is the dynamic viscosity of the gas mixture. 
Velocity inlet–pressure outlet boundary conditions were used for the system \Cref{eq:Stokes}. 

Thus, at the first step of the overall solution procedure,
the velocity field ${\bf v}$ was determined as a solution of system \Cref{eq:Stokes}. In the second step, due to (ii), using the predetermined velocity field, the concentration field $c$ and the integral \Cref{eq:bcurve} were computed. Both steps were performed with the PoreChem software \cite{PoreChem}. 

\section{Kernel Methods}\label{sect:kernels}

Kernel methods \cite{wendland2005scattered} comprise versatile tools for scattered data approximation,
revolving around the notion of a symmetric kernel $k: \Omega_d \times \Omega_d \rightarrow \R$.
An important type of kernels is given by strictly positive definite kernels, 
i.e.,\ kernels such that the so called kernel matrix $k(X_N, X_N) = (k(\fx_i, \fx_j))_{i,j=1}^N$ is positive definite for any choice of pairwise distinct points $X_N = \{\fx_i\}_{i=1}^N \subset \Omega_d$, $N \in \N$. In the context of machine learning the set $X_N$ is often referred to as the training set $\train{X}$, i.e., $X_N = \train{X}.$
In the following, we focus on the popular class of radial basis function kernels on $\Omega_d \subseteq \R^d$,
which can be expressed as $k(\fx, \fx') = \phi(\Vert \fx - \fx' \Vert)$ using a radial basis function $\phi: \R \rightarrow \R$.
Popular instances of such kernels are given by the following,
    \begin{align*}
    k(\fx,\fx') &= e^{-\epsilon^2||\fx - \fx'||^2} && \text{Gaussian kernel,} \\
    k(\fx,\fx') &= (1 + \epsilon  ||\fx - \fx'||) e^{-\epsilon||\fx - \fx'||} && \text{Matérn 1,} \\
    k(\fx,\fx') &= (3 + 3\epsilon||\fx - \fx'|| + \epsilon^2||\fx - \fx'||^2) e^{-\epsilon||\fx - \fx'||} && \text{Matérn 2,}
    \end{align*}
which are also later on used in \Cref{sect:num}.
The parameter $\epsilon$ is the so-called shape parameter, which allows to tune the width of these kernels.

Given not only input data $\{ \fx_i \}_{i=1}^N$, but also corresponding (possibly vector-valued) target values $\{ \ff_i \}_{i=1}^N \subset \R^b$, $b \geq 1$,
a well known representer theorem states that the optimal least squares kernel approximant is given by 
    \begin{align}
    \label{eq:kernel_modell_full}
        s_N(\fx) = \sum_{i=1}^N \alpha_i k(\fx, \fx_i).
    \end{align}
    with coefficients $\{ \alpha_j \}_{j=1}^N \subset \R^b$.
These coefficients can be computed directly by solving the regularized linear equation system $(k(X_N, X_N) + \lambda I)\alpha = y$,
where $\alpha \in \R^{N \times b}$ and $y \in \R^{N \times b}$ refer to a collection of the coefficients and target values in arrays. 
The regularization parameter $\lambda \geq 0$ allows to steer the robustness to outliers/noise versus the approximation accuracy in the training points. The value $\lambda=0$ corresponds to kernel interpolation. For this case $\lambda = 0$, 
the assumption on the strict positive definiteness of the kernel still ensures the solvability of the system.

\textbf{Greedy kernel approximation:}
In order to obtain a sparse kernel model, 
one strives to have a smaller expansion size $n \ll N$ within Eq.~\eqref{eq:kernel_modell_full}.
For this, greedy algorithms can be leveraged, which select a suitable subset of centers $X_n$ from $X_N$.
For this, they start with $X_0 = \{ \}$ and iteratively add points from $X_N$ to $X_n$ as $X_{n+1} := X_n \cup \{ \fx_{n+1} \}$ according to some selection criterion.
While there are several selection criterion in use in the literature,
we focus on the $f$-greedy selection criterion which reads
\begin{align*}
\fx_{n+1} := \mathrm{argmax}_{\fx_i \in X_N \setminus X_n} |\ff_i-s_n(\fx_i)|.
\end{align*}
This residual based selection criterion incorporates the data point of the largest error, 
i.e.\ directly aims at minimizing the maximal absolut error of the kernel model.
An efficient implementation of such greedy kernel algorithms is provided by the VKOGA package \cite{Santin2021VKOGA},
and a comprehensive analysis of the convergence of such greedy kernel algorithms was provided in \cite{wenzel2023analysis}.

\textbf{Two-layered kernels:}
In order to incorporate feature learning into kernel models, 
we make use of two-layered kernels \cite{wenzel2024data, wenzel2024application}:
These make use of a base kernel $k$ as given above, and combine it with a matrix $\fA \in \R^{d \times d}$,
such that the two-layered kernel is given by 
\begin{align}
    \label{eq:two_layered_kernel}
    k_\fA(\fx, \fx') = k(\fA\fx, \fA\fx') = \phi(\Vert \fA(\fx-\fx') \Vert).
\end{align}
With this, the two-layered kernel can be optimized to given data $(\{\fx_i\}_{i=1}^N, \{\ff_i\}_{i=1}^N)$ by optimizing the first layer matrix $\fA \in \R^{d \times d}$. 
Thus this two-layered structure may considerably improve the performance of the kernel model, 
especially for medium- to high-dimensional input data,
where an equal importance of all features is highly unlikely. The strength of two-layered kernel models and their superior performance over shallow kernel models has been observed, for example, when applied to heat conduction \cite{wenzel2024application} or when used as a surrogate model for chemical kinetics \cite{doeppel2024goal}. We will investigate whether this behaviour can also be observed for our current problem.
For the optimization of the matrix $\fA$, we employ the fast gradient descent based mini-batch optimization proposed in \cite{wenzel2024data} and extended it to vector valued target values $b > 1$ in \cite{wenzel2024application}.
The overall idea is to leverage an efficiently computable leave-one-out cross validation error loss,
which thus makes the kernel generalize well to unseen data.
In particular, this cross validation error loss makes use of both input and target values, i.e.,\ it is a supervised optimization.

The matrix $\fA$ makes the two-layered kernel $k_\fA$ even more interpretable:
The large singular values with corresponding right singular vectors within the singular value decomposition of the matrix $\fA$ correspond to the more important features within the data set,
while the smaller singular values with corresponding right singular vectors correspond to the less important features within the data set.

In the following \Cref{sect:num},
we make use of the notion ``single layered kernel'' to refer to standard kernel $k$,
while we make use of the notion ``two-layered kernel'' to refer to (optimized) kernels $k_\fA$.

\section{Numerical Experiments}\label{sect:num}
In this section we present the numerical experiments.
The goal is to approximate breakthrough curves \Cref{eq:bcurve} from voxel data which characterizes the geometry of the porous medium.
As explained in \Cref{sect:problem}, the voxel data is described by a vector $\fz \in \{0,1\}^{2 N_V}$.
In our numerical experiment, we choose $N_V = 150^3$.
We discretize a breakthrough curve $a(t)$ on an equidistant temporal grid $t_i$, $i=1,\dots,\numTimesteps$, as the vector $\fa \vcentcolon= (a(t_i))_{i=1}^{\numTimesteps} \in \R^{\numTimesteps}$ with $\numTimesteps=500$ time steps.

\begin{figure}[t]
    \centering
    \graphicspath{{Plots/models}}
    \def\svgwidth{.9\textwidth}
    
\begingroup%
  \makeatletter%
  \providecommand\color[2][]{%
    \errmessage{(Inkscape) Color is used for the text in Inkscape, but the package 'color.sty' is not loaded}%
    \renewcommand\color[2][]{}%
  }%
  \providecommand\transparent[1]{%
    \errmessage{(Inkscape) Transparency is used (non-zero) for the text in Inkscape, but the package 'transparent.sty' is not loaded}%
    \renewcommand\transparent[1]{}%
  }%
  \providecommand\rotatebox[2]{#2}%
  \newcommand*\fsize{\dimexpr\f@size pt\relax}%
  \newcommand*\lineheight[1]{\fontsize{\fsize}{#1\fsize}\selectfont}%
  \ifx\svgwidth\undefined%
    \setlength{\unitlength}{413.08295303bp}%
    \ifx\svgscale\undefined%
      \relax%
    \else%
      \setlength{\unitlength}{\unitlength * \real{\svgscale}}%
    \fi%
  \else%
    \setlength{\unitlength}{\svgwidth}%
  \fi%
  \global\let\svgwidth\undefined%
  \global\let\svgscale\undefined%
  \makeatother%
  \begin{picture}(1,0.58685955)%
    \lineheight{1}%
    \setlength\tabcolsep{0pt}%
    \put(0.11239177,0.56320329){\makebox(0,0)[t]{\lineheight{0.94999999}\smash{\begin{tabular}[t]{c}\textbf{input}\end{tabular}}}}%
    \put(0.88655082,0.56320329){\makebox(0,0)[t]{\lineheight{0.94999999}\smash{\begin{tabular}[t]{c}\textbf{output}\end{tabular}}}}%
    \put(0,0){\includegraphics[width=\unitlength,page=1]{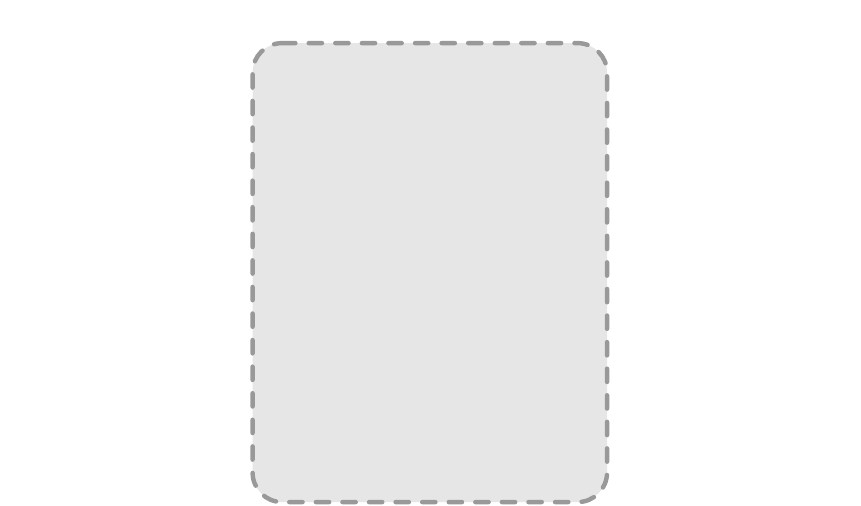}}%
    \put(0.49932639,0.56320329){\makebox(0,0)[t]{\lineheight{0.94999999}\smash{\begin{tabular}[t]{c}\textbf{model}\end{tabular}}}}%
    \put(0.6027782,0.31840585){\makebox(0,0)[t]{\smash{\begin{tabular}[t]{c}kernel\\model\end{tabular}}}}%
    \put(0.6027782,0.05306152){\makebox(0,0)[t]{\smash{\begin{tabular}[t]{c}kernel\\model\end{tabular}}}}%
    \put(0,0){\includegraphics[width=\unitlength,page=2]{Plots/models/num_exp_approaches_structued.pdf}}%
    \put(0.11255603,0.50733704){\makebox(0,0)[t]{\smash{\begin{tabular}[t]{c}voxel data\\($2\cdot 150^3 \times 1$)\end{tabular}}}}%
    \put(0,0){\includegraphics[width=\unitlength,page=3]{Plots/models/num_exp_approaches_structued.pdf}}%
    \put(0.41083946,0.32237116){\makebox(0,0)[t]{\smash{\begin{tabular}[t]{c}extract\\morphological\\features\end{tabular}}}}%
    \put(0,0){\includegraphics[width=\unitlength,page=4]{Plots/models/num_exp_approaches_structued.pdf}}%
    \put(0.41144101,0.05141294){\makebox(0,0)[t]{\smash{\begin{tabular}[t]{c}extract\\PCA features\end{tabular}}}}%
    \put(0,0){\includegraphics[width=\unitlength,page=5]{Plots/models/num_exp_approaches_structued.pdf}}%
    \put(0.2981774,0.36720945){\makebox(0,0)[t]{\lineheight{1.25}\smash{\begin{tabular}[t]{c}\textbf{(MF) }\end{tabular}}}}%
    \put(0,0){\includegraphics[width=\unitlength,page=6]{Plots/models/num_exp_approaches_structued.pdf}}%
    \put(0.2981774,0.09850586){\makebox(0,0)[t]{\lineheight{1.25}\smash{\begin{tabular}[t]{c}\textbf{(PCA) }\end{tabular}}}}%
    \put(0.88690631,0.50733704){\makebox(0,0)[t]{\smash{\begin{tabular}[t]{c}breakthrough\\curve\\($500 \times 1$)\end{tabular}}}}%
  \end{picture}%
\endgroup%
    \caption{Two feature extraction strategies based models to approximate breakthrough curves from voxel data: Model with morphological features (MF), and model with PCA features (PCA).}
    \label{fig:methods}
\end{figure}

We compare two different kernel-based models to learn the breakthrough curves (see \Cref{fig:methods}): (MF) applies morphological features prescribed by expert knowledge (see \cref{sect:appendix_S}), and 
(PCA) is based on features that are learned with the principal component analysis. 
The kernel function $k$ is chosen to be either a shallow one-layered or a two-layered kernel.
We refer to these models in the scope of our paper as MF-1L-kernel and MF-2L-kernel or PCA-1L-kernel and PCA-2L-kernel depending on the depth of the kernel (one-layered/two-layered) and whether we extract morphological features (MF) or PCA features.

Both types of models (MF) and (PCA) are based on a feature map $\featureMap: \R^{\numGeom} \to \R^{\numFeatures}$ with $\numFeatures$ features and a kernel function $k: \R^{d} \times \R^{d} \to \R$ with $d = \numFeatures$,
such that the resulting kernel models using the single-layer or the two-layered kernels are given as
\begin{align}
    s_n(\fz) = \sum\limits_{i=1}^n\alpha_i k(\featureMap(\fz), \featureMap(\fz_i)), \qquad \text{respectively} \qquad 
    s_n(\fz) = \sum\limits_{i=1}^n\alpha_i k(\fA \featureMap(\fz), \fA \featureMap(\fz_i)),\label{eqn:feature_model}
\end{align}
with coefficients $\alpha_i \in \R^{b}$, $b = \numTimesteps$, and centers $\fz_i \in \{0,1\}^{\numGeom}$ for $i=1,\dots,n$. The expansion size $n$ is fixed to $n = 10.$ Choosing smaller values for $n$ would worsen the results presented in the next subsections considerably while increasing $n$ would only influence the approximation quality on the test set slightly.

All of these models involve hyper-parameters, which are listed in \Cref{tab:cv}. 
Suitable values for the hyper-parameters are determined via a leave-one-out cross-validation (LOOCV) on the training data set. Note that in the cases where two-layered kernels are applied, no LOOCV for the shape parameter is performed and instead the matrix $\fA$ is optimized using the procedure described in the previous section. The kernel function and the regularization parameter are determined via LOOCV in both cases.
\begin{table}[t]
    \centering
    \begin{tabular}{ll}
    \toprule
         \textbf{hyper-parameter}& \textbf{values} \\ \midrule
         kernel fun.\ & Matérn 1, Matérn 2, Gauß\\\midrule
         Shape parameters& 1, 1/2, 1/4, 1/8, 1/16, 1/32\\\midrule
         Regularization parameters& 0, $10^{-2}, 10^{-3}, 10^{-4}, 10^{-5}, 10^{-6}$ \\ \bottomrule
    \end{tabular}
    \caption{Hyper-parameters of the kernel function $k$, shape parameter $\epsilon$ and regularization parameter $\lambda$ used during LOOCV of the MF-based models and PCA-based models.}
    \label{tab:cv}
\end{table}

In the following,
we discuss the generation of the training and test data (\Cref{sect:data}),
the training and results of the (MF) models based on morphological features (\Cref{sect:mf}), and of the (PCA) models based on PCA features (\Cref{sect:pca}).

\subsection{Training and Test Data}\label{sect:data}
\begin{figure}[t]
    \centering
    \begin{minipage}{0.6\textwidth}
        \includegraphics[width=\textwidth]{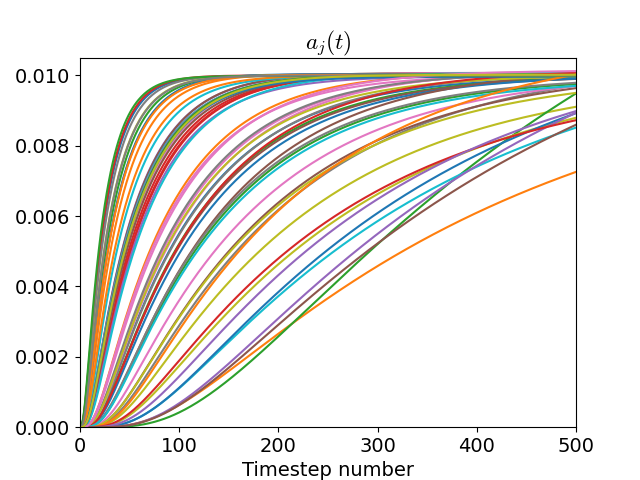}
    \end{minipage}
    \caption{All breakthrough curves $a_i(t)$ computed from the voxel data $\fz_i$ for $i\in \{1, \dots, 59 \}$.} 
    \label{fig:allBCurves}
\end{figure}

In total, we consider $\numSamples=59$ samples $X\vcentcolon=\{\fz_i\}_{i=1}^{\numSamples}$ of voxel data $\fz_i$ and the corresponding time-discretized breakthrough curves $\fa_i$ (see \Cref{fig:allBCurves}).
The breakthrough curves $\fa_i$ are obtained by solving \Cref{eq:CDR} and \Cref{eq:Stokes} based on the geometry described by the corresponding voxel data $\fz_i$.
The average computational cost of solving \Cref{eq:CDR} and \Cref{eq:Stokes} for one sample is $7514$ wall-clock seconds (wcs) with a standard deviation of $286$ wcs on a workstation with two Intel Xeon Gold 6240R (48 cores in total).

We use an approximate $80\%$--$20\%$ training--test split, which results in $ \numSamplesTrain = N = 47$ training samples $\train{X} = X_N$ and $\numSamplesTest = 12$ test samples $\test{X}$ with $X = \train{X} \cup \test{X}$ and $\train{X} \cap \test{X} = \emptyset$. 
To demonstrate robustness with respect to the choice of the training--test split, we consider three different random splits, where we use the same three random splits to measure the performance of each model.
As an error measure we use a relative error on the test set
\begin{equation}
    \erel = \frac{1}{|{\test{X}}|}\sqrt{\sum\limits_{\fz_i \in {\test{X}}}\frac{||\fa_i - s_N(\fz_i)||^2}{||\fa_i||^2}}. \label{eq:e_rel}
\end{equation}

\subsection{Kernel Models on Morphological Features}\label{sect:mf}
Because of the high dimensionality of the input space, it is necessary to apply some reduction technique before training the machine learning model.
In this section we extract $n_f = 6$  morphological geometry features that describe the porous medium. 
These are the porosity, the washcoat volume fraction, the volume of free pores, the surface area for the free pores, the integral of mean curvature of the free pores and the integral of total curvature of the free pores. 
For further information on how they are computed see supplementary \Cref{sect:appendix_S}. 
Thus, the machine learning model can be represented by setting $\Phi = \Phi_\text{MF}$ in \Cref{{eqn:feature_model}}
$\Phi_\text{MF}:\R^{2N_V}\rightarrow\R^{n_f}$ the mapping that computes the morphological features from a given geometry. 

\begin{figure}[t]
    \centering
    \begin{subfigure}{0.475\textwidth}
\includegraphics[width=\textwidth]{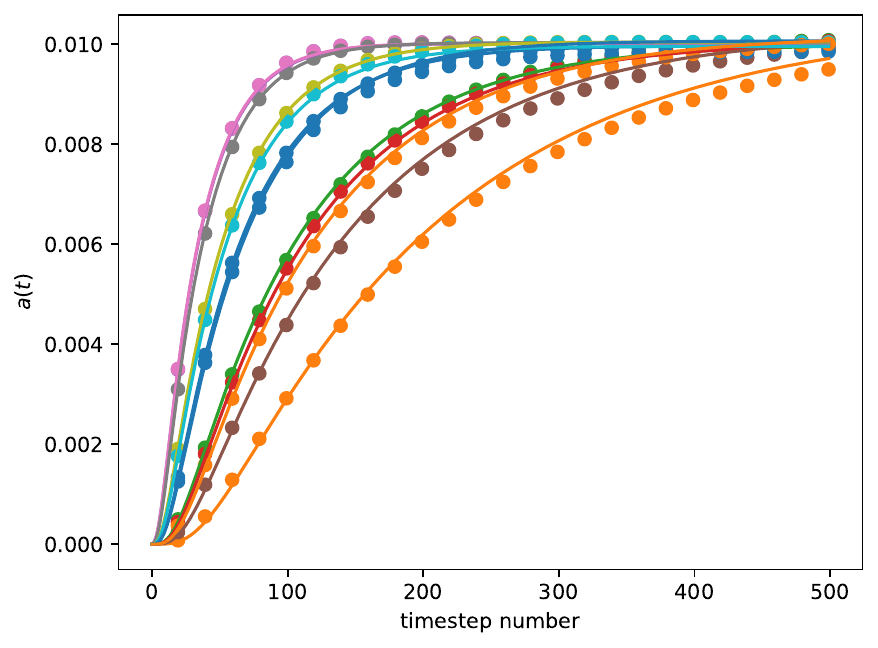}
\caption{Split 0}
\label{fig:MF_1La}
\vspace{1em}
\end{subfigure}
\begin{subfigure}{0.475\textwidth}
\includegraphics[width=\textwidth]{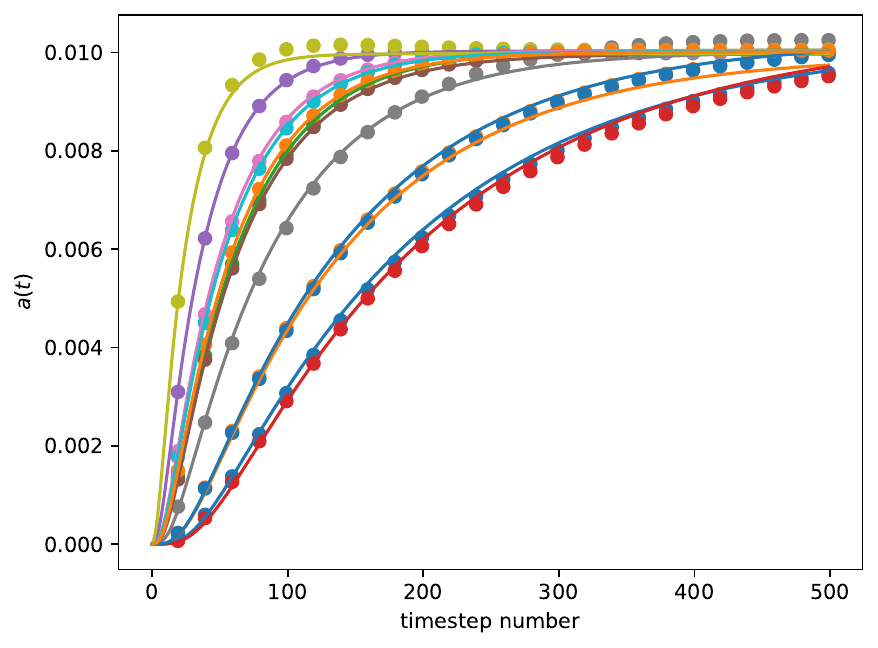}
\caption{Split 1}
\label{fig:MF_1Lb}
\vspace{1em}
\end{subfigure}

\begin{subfigure}{0.475\textwidth}
\includegraphics[width=\textwidth]{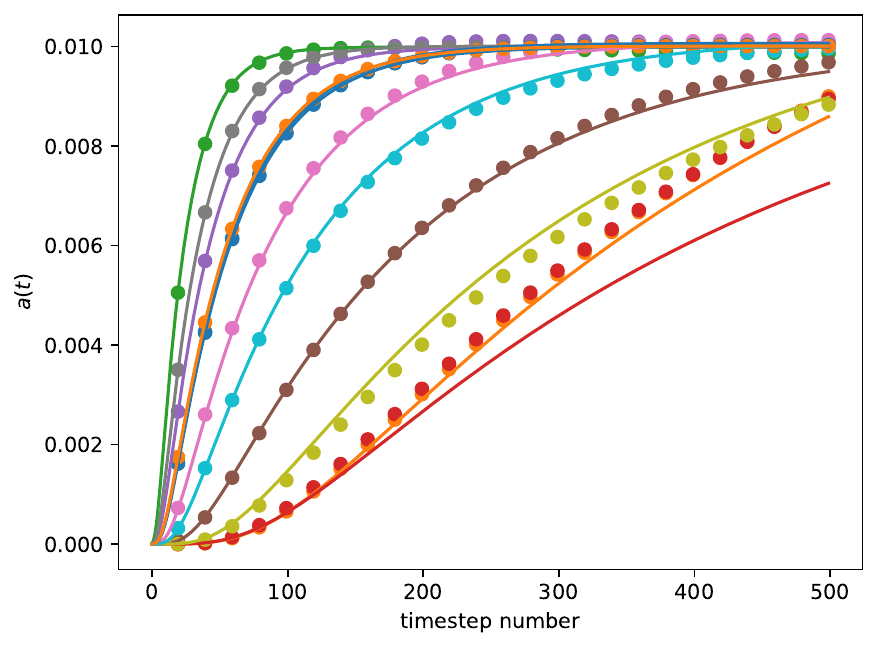}
\caption{Split 2}
\vspace{1em}
\label{fig:MF_1Lc}
\end{subfigure}
\caption{MF-1L-kernel: Predicted breakthrough curves on the test set $\test{X}$ for three different training--test splits. Solid: Breaktrough curves, dotted: predictions. The different colours correspond to the different breakthrough curves.}
\label{fig:MF_1L}
\end{figure}

For the first experiment we use shallow kernels and present the approximated breakthrough curves in \Cref{fig:MF_1L}. 
We observe that, except for the red outlier curve in the bottom diagram, all curves are well approximated. This corresponds to the relative test error presented in 
\Cref{tab:MF_1L}, where we see that for the first two data splits a relative error of about $0.01\%$ can be achieved whereas, for the third split, we only achieve an error of about $0.42\%$. We further observe from \Cref{tab:MF_1L} that for the first two data splits exactly the same hyper-parameters (kernel function, shape parameter and regularization parameter) are chosen. However, they are chosen differently for the third data set, which may be due to the red outlier curve being part of the test set for the third data split.

\begin{table}[t]
    \centering
    \begin{tabular}{ccccc}
\toprule
split & rel.~err. & kernel fun.\ & shape parameter & reg.\ par.\ \\
\midrule
0 & 1.49$\cdot10^{-3}$ & Matérn 1 & 1.00 & 1.00$\cdot10^{-2}$ \\
1 & 1.60$\cdot10^{-4}$ & Matérn 1 & 1.00 & 1.00$\cdot10^{-2}$ \\
2 & 4.23$\cdot10^{-2}$ & Matérn 1 & 2.50$\cdot10^{-1}$ & 1.00$\cdot10^{-4}$ \\
\bottomrule
\end{tabular}
 \caption{MF-1L-kernel: Relative error \eqref{eq:e_rel} and the selected hyper-parameters from \Cref{tab:cv} for three different training--test splits (0-2).}
    \label{tab:MF_1L}
\end{table}

\begin{figure}[t]
    \centering
    \begin{subfigure}{0.475\textwidth}
\includegraphics[width=\textwidth]{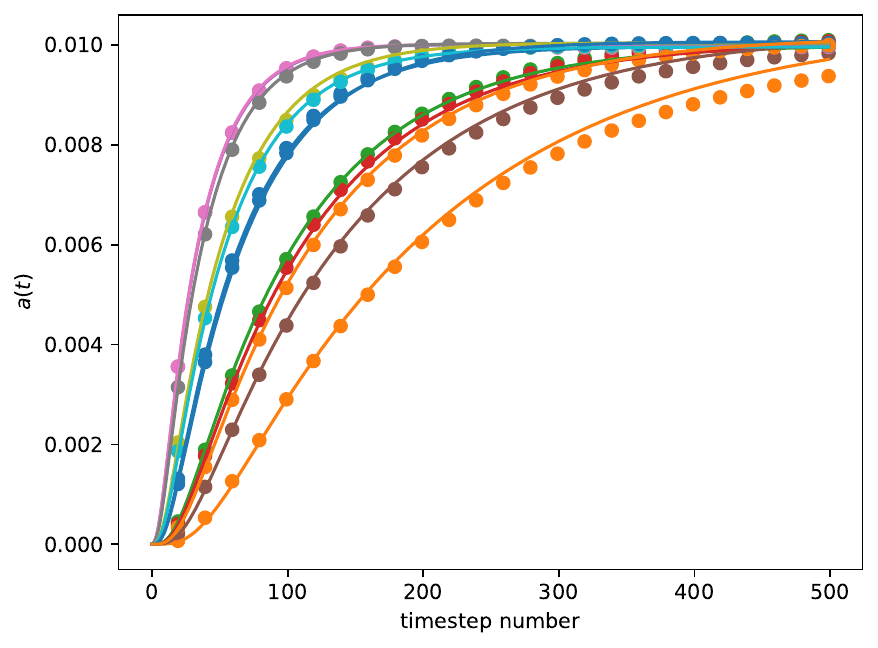}
\label{fig:MF_2La}
\caption{Split 0}
\vspace{1em}
\end{subfigure}
\begin{subfigure}{0.475\textwidth}
\includegraphics[width=\textwidth]{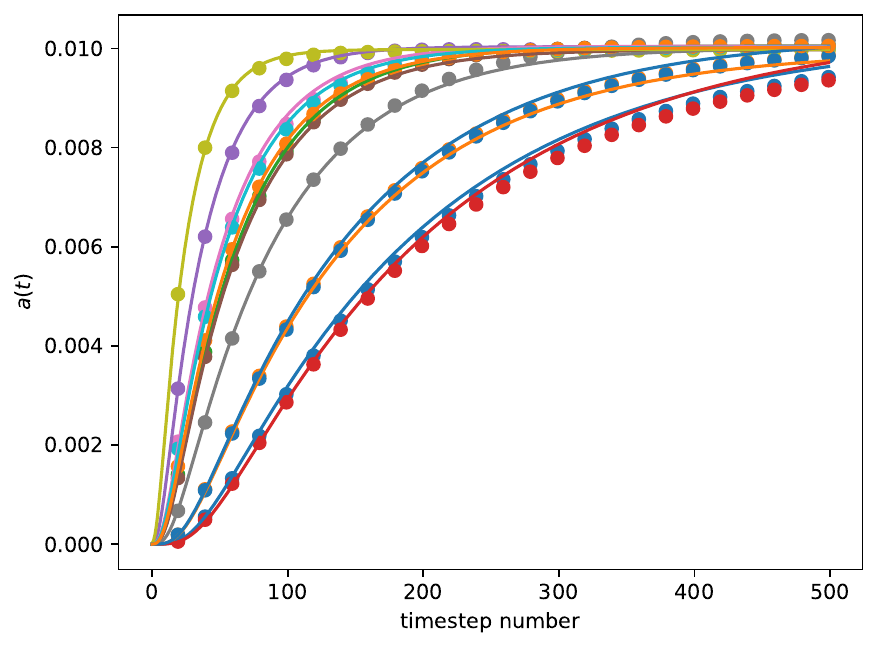}
\label{fig:MF_2Lb}
\caption{Split 1}
\vspace{1em}
\end{subfigure}
\begin{subfigure}{0.475\textwidth}
\includegraphics[width=\textwidth]{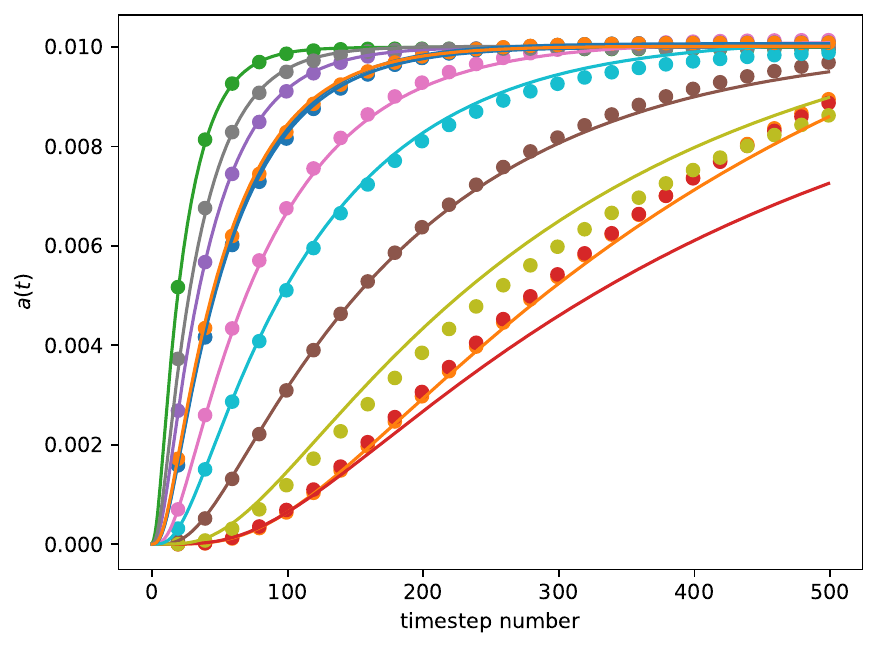}
\label{fig:MF_2Lc}
\caption{Split 2}
\vspace{1em}
\end{subfigure}
\caption{MF-2L-kernel: Predicted breakthrough curves on the test set $\test{X}$ for three different training--test splits. Solid: Breaktrough curves, dotted: predictions. The different colours correspond to the different breakthrough curves.}
\label{fig:MF_2L}
\end{figure}

For the next experiment we included two-layer optimization of the kernel model and present the approximated breakthrough curves in \Cref{fig:MF_2L}.
We observe that a similar approximation quality can be achieved using two-layered kernels. All breakthrough curves except for the red outlier curve in the bottom diagram are well approximated. The relative test error (see \Cref{tab:MF_2L}) is slightly improved for the third data split. However, for the the first and second split, it gets slightly worse. We observe that similar to the experiment with shallow kernel models, for each data split, the Matérn 1 kernel is selected and that, for the third split, a smaller regularization parameter is selected than for the first and second one. To compare the shape transformations by the first layer, we present the first layer matrix $\fA$ in \Cref{tab:MF_2L}. For an easier comparison the matrix $\fA$ is visualized using the color-map at the bottom of the table.
We observe that the matrices look very similar for each data split. 

\setlength{\tabcolsep}{4pt}
\setlength{\arraycolsep}{1pt}
\begin{table}[ht!]
    \centering
    \begin{tabular}{ccccc}
\toprule
split & rel.~err. & kernel fun.\ & $\fA$ & reg.\ par.\ \\
\midrule 
0 & 1.54$\cdot10^{-4}$ & Matérn 2 & \begin{minipage}{.2\textwidth}\includegraphics[scale=0.17]{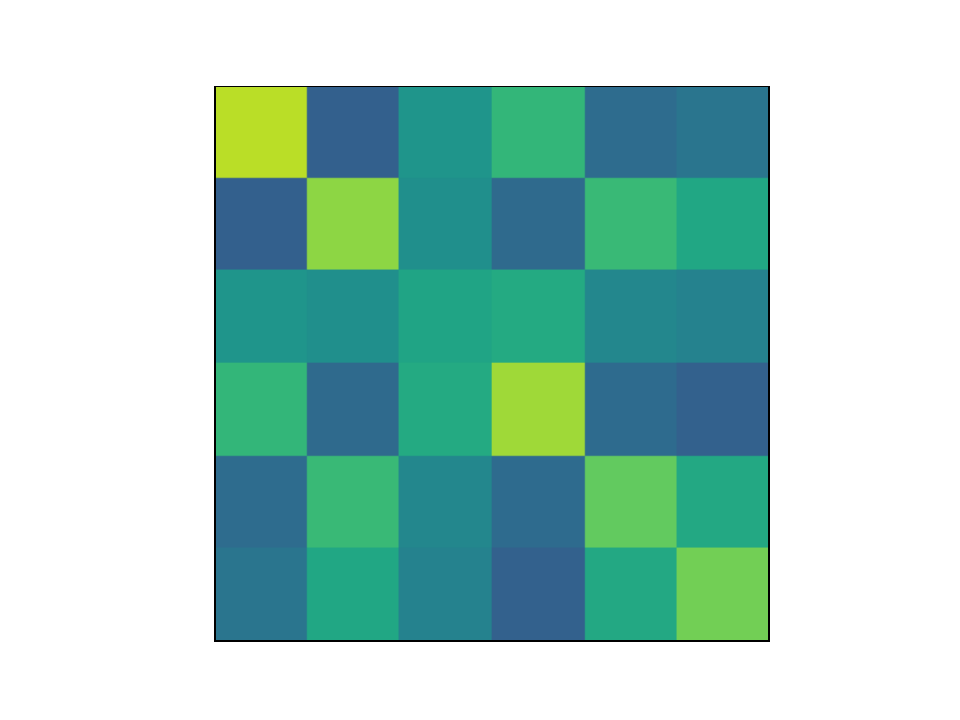} \end{minipage}
& 1.00$\cdot10^{-3}$ \\
1 & 2.37$\cdot10^{-4}$ & Matérn 2&  
\begin{minipage}{.2\textwidth}\includegraphics[scale=0.17]{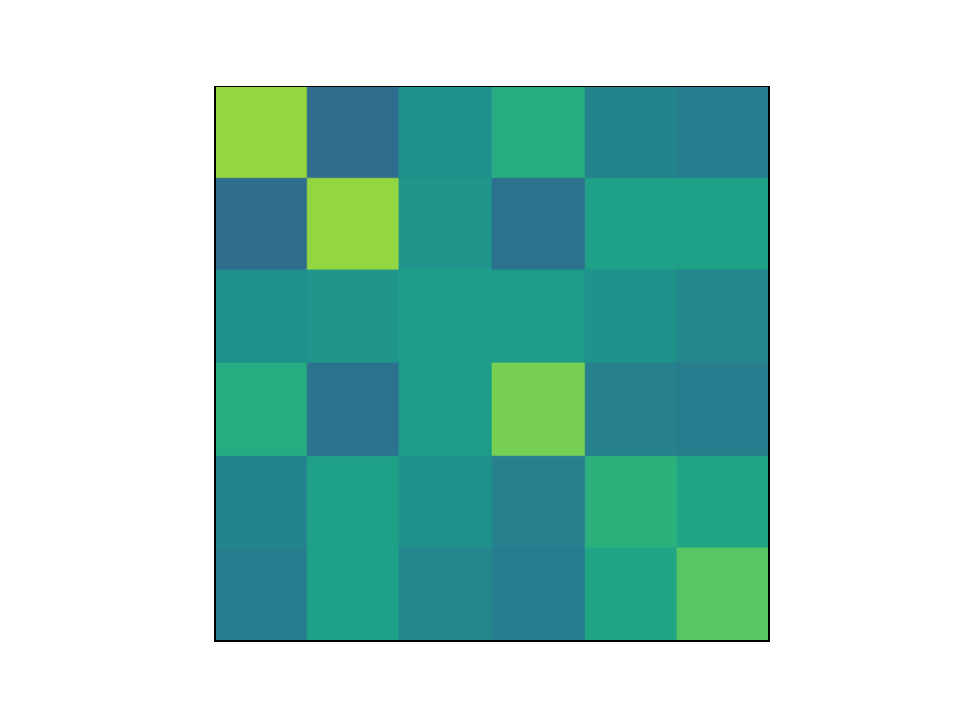} \end{minipage}
 & 1.00$\cdot10^{-3}$ \\ 
2 & 3.99$\cdot10^{-3}$ & Matérn 1 & \begin{minipage}{.2\textwidth}\includegraphics[scale=0.17]{Plots/MF2.pdf} \end{minipage}
 & 1.00$\cdot10^{-4}$ \\
\bottomrule \\[-1em]
\multicolumn{5}{r}{%
\includegraphics[scale = 0.15]{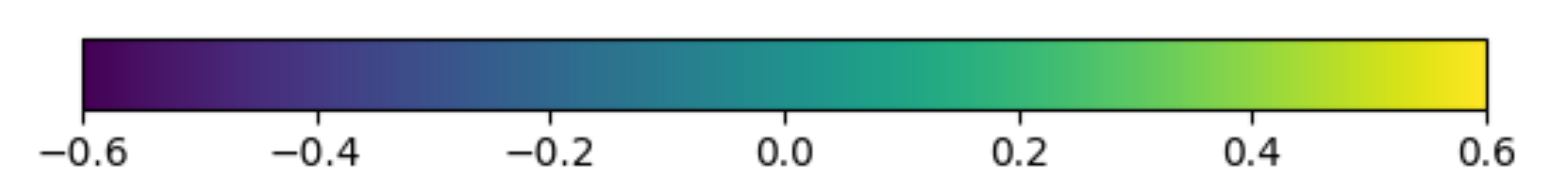}}
\end{tabular}
    \caption{MF-2L-kernel: Relative error \eqref{eq:e_rel} and the selected hyper-parameters from \Cref{tab:cv} for three different training--test splits (0-2).}
    \label{tab:MF_2L}
\end{table}

\subsection{Kernel Models on PCA Features}\label{sect:pca}
Next, we use the PCA to define a feature map $\featureMap$.
In order to be comparable to the previous experiment, we choose the number of features as $\numFeatures = 6$.
Following the idea of the PCA \cite{J02}, we choose the PCA feature map $\featureMapPCA(\fz) = \rT\fU_{\numFeatures} \fz$ based on the first $\numFeatures$ left-singular vectors $\fU_{\numFeatures} \in \R^{\numGeom \times \numFeatures}$ of the matrix\footnote{Technical note: The voxel data is saved as a boolean array $\fz\in \{0,1\}^{\numGeom}$. To compute the SVD, we convert this data to a floating point number which is why $\fZ \in \R^{\numGeom \times \numSamplesTrain}$.} $\fZ \vcentcolon= (\fz)_{\fz \in \train{X}} \in \R^{\numGeom \times \numSamplesTrain}$ via the SVD
\begin{align}
    &\fZ = \fU \fSigma \rT\fV&
&\text{with}&
&\fU_{\numFeatures} \vcentcolon= \fU(:, :\numFeatures)
\end{align}
and set $\Phi = \featureMapPCA$ in \Cref{eqn:feature_model}.

\begin{figure}[t]
\centering
\begin{subfigure}{0.475\textwidth}
\includegraphics[width=\textwidth]{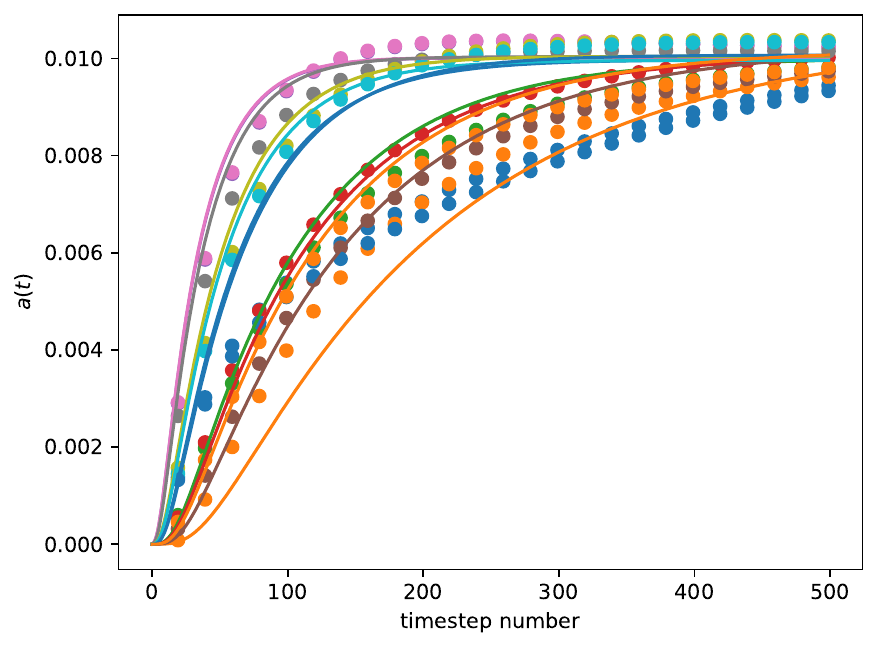}
\label{fig:PCA_1La}
\caption{Split 0}
\vspace{1em}
\end{subfigure}
\begin{subfigure}{0.475\textwidth}
\includegraphics[width=\textwidth]{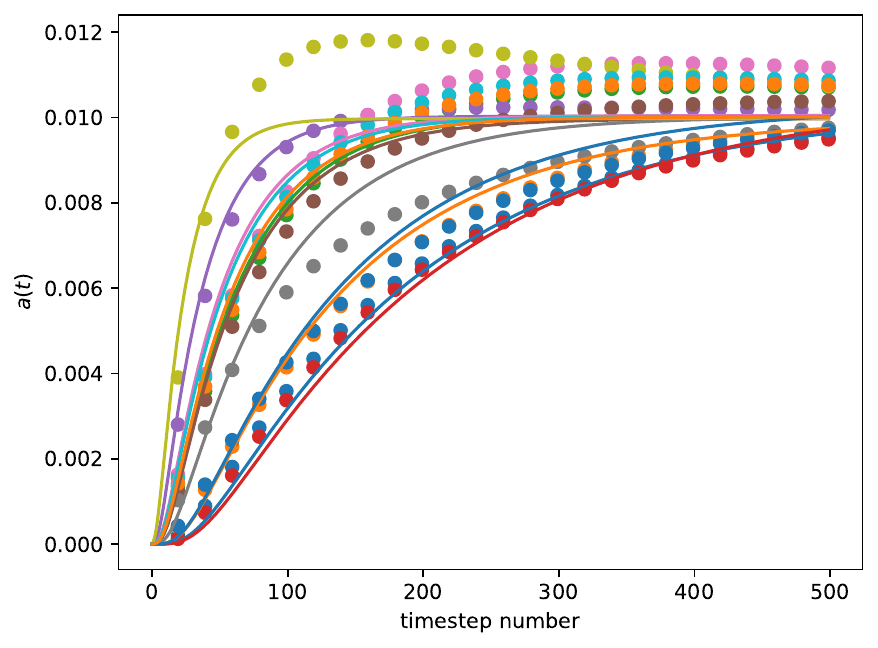}
\label{fig:PCA_1Lb}
\caption{Split 1}
\vspace{1em}
\end{subfigure}
\begin{subfigure}{0.475\textwidth}
\includegraphics[width=\textwidth]{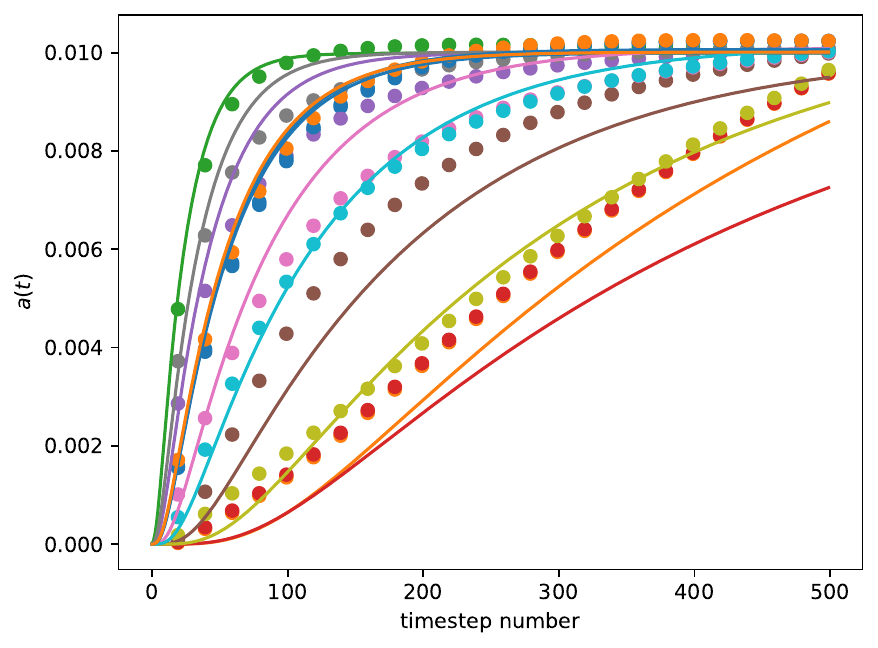}
\label{fig_PCA_1Lc}
\caption{Split 2}
\vspace{1em}
\end{subfigure}
\begin{subfigure}{0.475\textwidth}
\includegraphics[width = \textwidth]{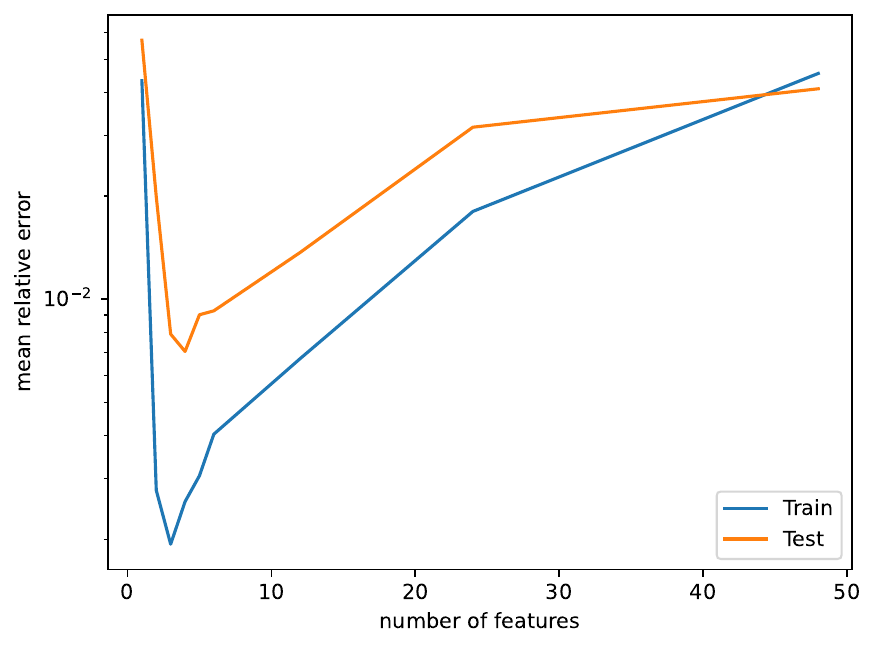}
\label{fig:PCA_1Ld}
\caption{Train and Test error over $\numFeatures$}
\vspace{1em}
\end{subfigure}
\caption{PCA-1L-kernel: \textbf{(a)-(c)}: Predicted breakthrough curves on the test set $\test{X}$ for $\numFeatures = 6$ and three different training--test splits. Solid: Breaktrough curves, dotted: predictions. The different colours correspond to the different breakthrough curves. \textbf{(d)}: Relative test and train error (\Cref{{eq:e_rel}}) over $\numFeatures$.
}
\label{fig:PCA_1L}
\end{figure}
\begin{table}[t]
    \centering
    \begin{tabular}{ccccc}
\toprule
split & rel.~err. & kernel fun.\ & shape parameter & reg.\ par.\ \\
\midrule
0 & 9.62$\cdot10^{-3}$ & Matérn 1 & 6.25$\cdot10^{-2}$ & 1.00$\cdot10^{-3}$ \\
1 & 4.77$\cdot10^{-3}$ & Gauss & 6.25$\cdot10^{-2}$ & 1.00$\cdot10^{-4}$ \\
2 & 1.34$\cdot10^{-2}$ & Matérn 2 & 2.5$\cdot10^{-1}$ & 1.00$\cdot10^{-2}$ \\
\bottomrule
\end{tabular}

  \caption{PCA-1L-kernel: Relative error \eqref{eq:e_rel} and the selected hyper-parameters from \Cref{tab:cv} for three different training--test splits (0-2).}
    \label{tab:PCA_1L}
\end{table}

In \Cref{fig:PCA_1L} we present the prediction of the breakthrough curves for the PCA-1L-kernel model on the test set $\test{X}$ for the three different randomized training--test splits.
We observe that many of the breakthrough curves are well approximated and a mean relative error of about $0.5\% - 1.3\%$ (see \Cref{tab:PCA_1L}) can be achieved. However, some breakthrough curves are badly approximated, especially the two blue curves for the first data split. We observe that three different kernels, three different shape parameters and three different regularisation parameters are chosen for each data split. 
We further observe from \Cref{fig:PCA_1Ld} $\numFeatures = 6$ is also a suitable choice in the sense that a small test terror is achieved for that setting compared to higher and smaller values of $\numFeatures$.

\begin{figure}[t]
    \centering
    \begin{subfigure}{0.475\textwidth}
\includegraphics[width=\textwidth]{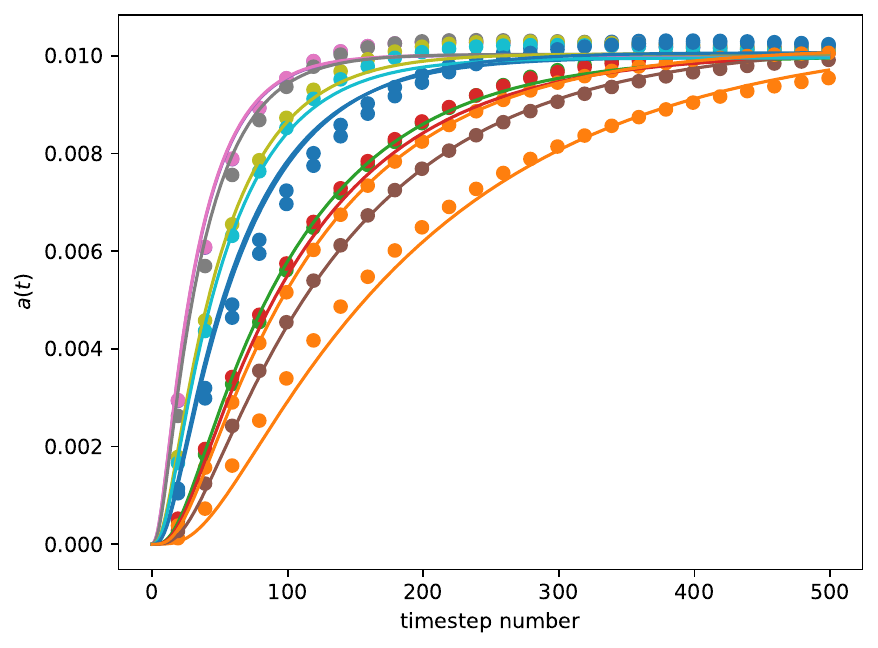}
\label{fif:PCA_2La}
\caption{Split 0}
\vspace{1em}
\end{subfigure}
\begin{subfigure}{0.475\textwidth}
\includegraphics[width=\textwidth]{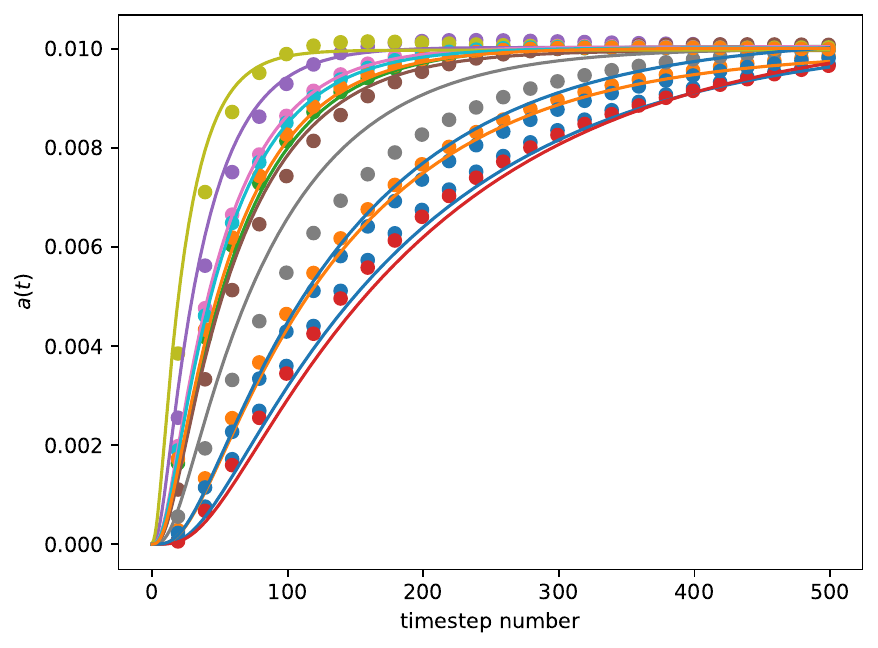}
\label{fig:PCA_2Lb}
\caption{Split 1}
\vspace{1em}
\end{subfigure}
\begin{subfigure}{0.475\textwidth}
\includegraphics[width=\textwidth]{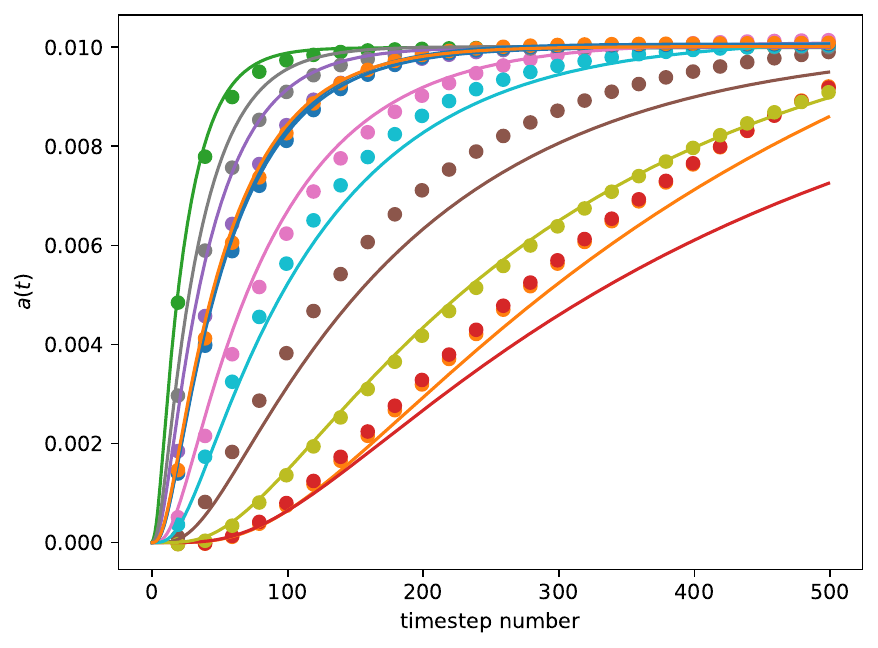}
\label{fig:PCA_2Lc}
\caption{Split 2}
\vspace{1em}
\end{subfigure}
\begin{subfigure}{0.475\textwidth}
\includegraphics[width = \textwidth]{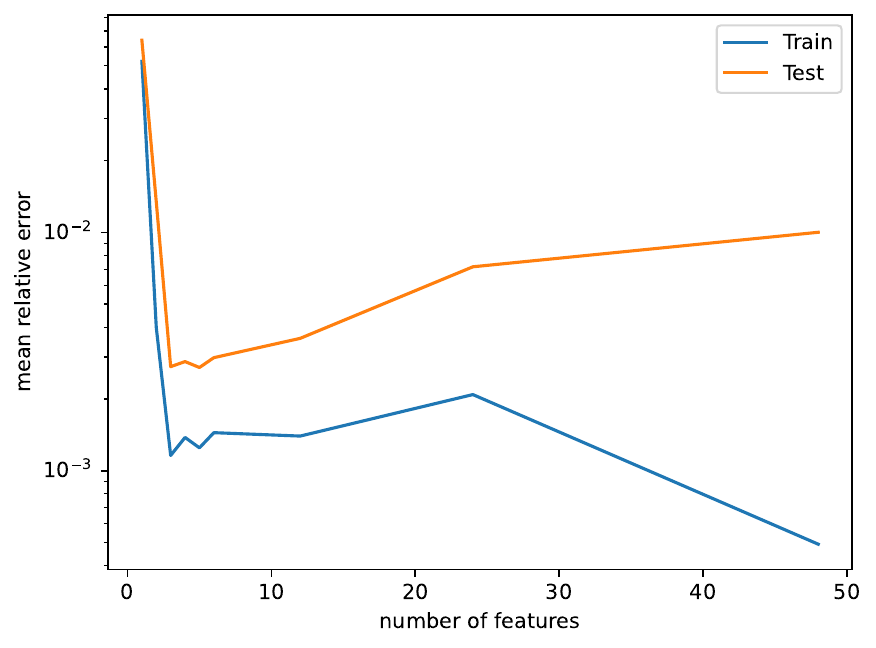}
\label{fig:PCA_2Ld}
\caption{Relative test and train error over $\numFeatures$}
\vspace{1em}
\end{subfigure}
\caption{\textbf{(a)-(c)}: Top \& left: Predicted breakthrough curves on the test set $\test{X}$ for $\numFeatures = 6$ and three different training--test splits. Solid: Breaktrough curves, dotted: predictions. The different colours correspond to the different breakthrough curves. \textbf{(d)}: Relative test and train error (\Cref{{eq:e_rel}}) over $\numFeatures$.}
\label{fig:PCA_2L}
\end{figure}

For the next experiment, we consider the PCA-2L-kernel model. The results for the analogous experiments to the previous paragraph are presented in \Cref{fig:PCA_2L}. 
We chose again $\numFeatures = 6$ to be consistent with the previous experiments. However, as we observe in \Cref{fig:PCA_2Ld} that for $\numFeatures = 3$ slightly better results could have been achieved. 

\begin{table}[ht!]
    \centering
    \begin{tabular}{ccccc}
\toprule

split & rel.~err. & kernel fun.\ & $\fA$ & reg.\ par.\ \\
\midrule
0 & 7.94$\cdot10^{-4}$ & Matérn 1 & \begin{minipage}{.2\textwidth}\includegraphics[scale=0.17]{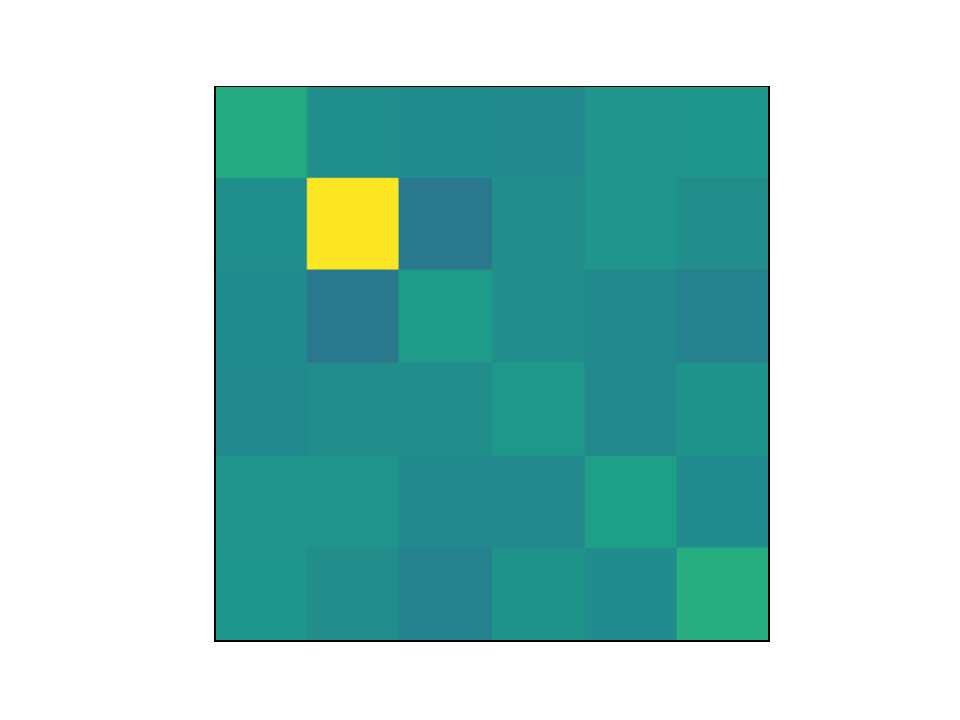} \end{minipage}
& 1.00$\cdot10^{-3}$ \\
1 & 1.10$\cdot10^{-3}$ & Matérn 2 &  \begin{minipage}{.2\textwidth}\includegraphics[scale=0.17]{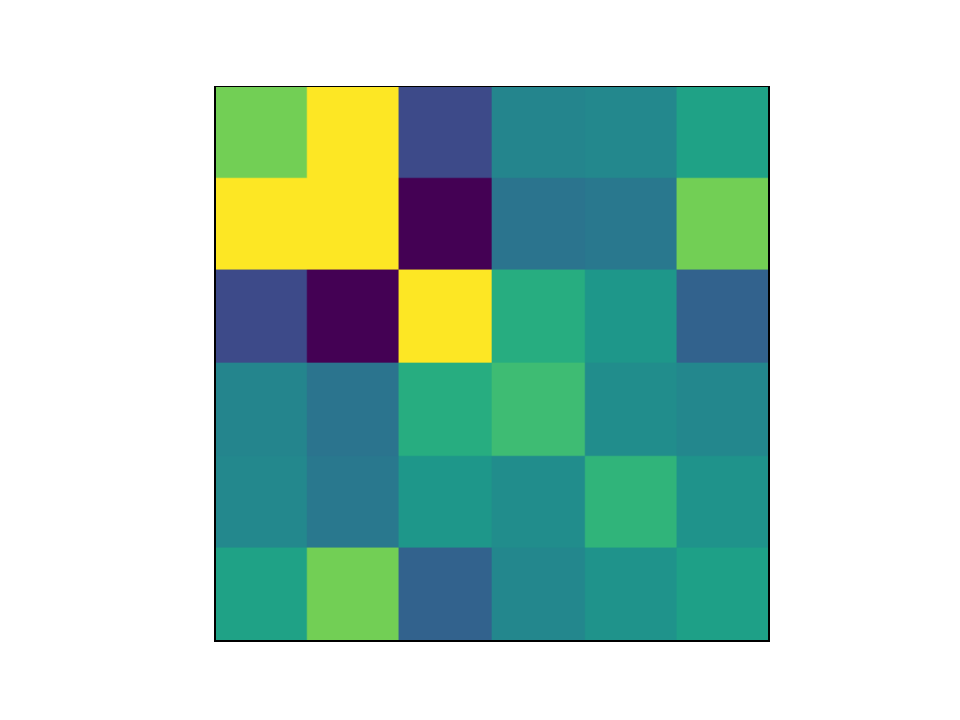} \end{minipage}
& 1.00$\cdot10^{-2}$ \\
2 & 7.05$\cdot10^{-3}$ & Matérn 1 & \begin{minipage}{.2\textwidth}\includegraphics[scale=0.17]{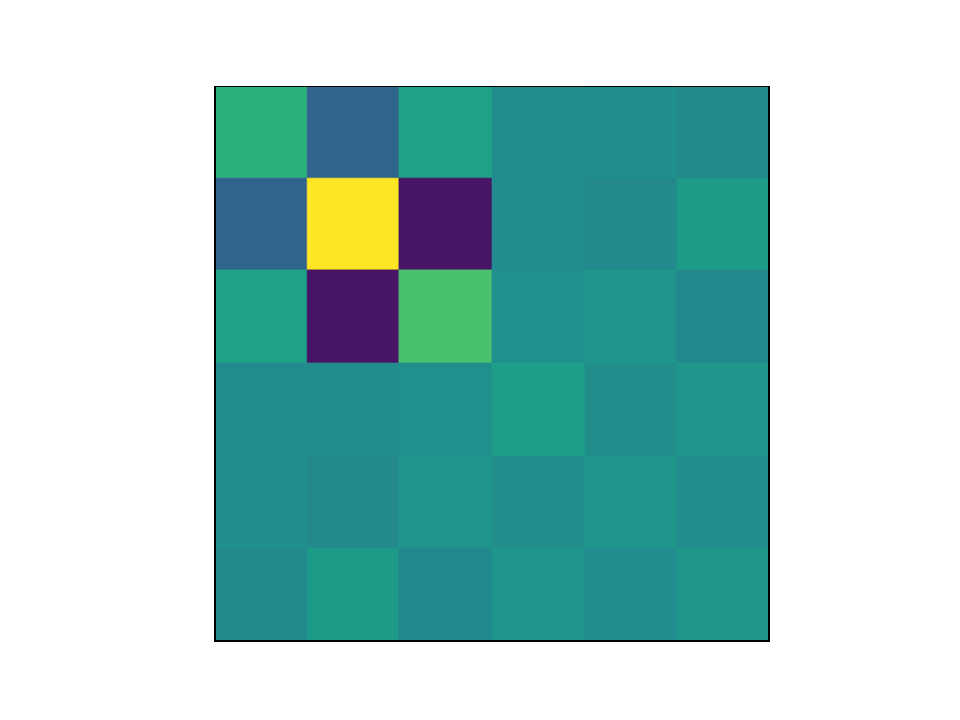} \end{minipage}
& 1.00$\cdot10^{-4}$ \\
\bottomrule \\[-1em]
\multicolumn{5}{r}{
\includegraphics[scale = 0.15]{Plots/colormap.png}}
\end{tabular}
    \caption{PCA-2L-kernel: Relative error \eqref{eq:e_rel} and the selected hyper-parameters from \Cref{tab:cv} for three different training--test splits (0-2).}
    \label{tab:PCA_2L}
\end{table}
We observe that the results from the shallow kernel models are considerably improved by using two-layered kernels. 
The relative errors from \Cref{tab:PCA_2L} can be reduced to $0.08\%-0.7\%$ and we do not observe badly approximated curves for the first data split anymore. 
Compared to the the MF-2L-kernel experiments from the previous section, we also observe that the approximation qualities are quite similar. 
This means that the PCA features transformed by the first layer of the two-layered kernel are able to describe the geometry almost as well as the morphological features.
This shows the strength of two-layered kernels, as the PCA feature extraction is purely data-based and no expert knowledge had to be used. 

We further observe from \Cref{tab:PCA_2L} that again similar hyper-parameters are selected for each data split. For example, the same kernel (Matérn 1) has been selected for the first and third data split. Moreover, we see similarities in the first layer matrices $\fA$ of the kernels. For all three matrices the entries in the left upper $3\times3$ block are considerably larger than the other entries which means that the first three modes are the important ones. Compared to \Cref{tab:MF_2L} there are larger differences between the first layers. This can be explained by the fact that, in contrast to the MF-feature map, the data-based PCA-feature map is different for all three data splits. This is due to the linear mappings defined by the PCA modes being different in each split. 

Lastly, we compare the mean errors averaged over the three different training--test splits in \Cref{tab:errors}. We observe, that the models based on morphological feature extraction work very well in combination with both one-layered and two-layered kernels. In contrast, combining one-layered kernels with PCA feature extraction yields an average error that is almost one order of magnitude larger. This error is considerably improved by applying two-layered kernels and almost the same accuracy as with the morphological feature extraction is achieved.

\begin{table}[t]
    \centering
    \begin{tabular}{c|c|c}
    \toprule
         & 1L & 2L\\ \midrule
       MF & 1.51$\cdot10^{-3}$ & 1.46$\cdot10^{-3}$ \\ \midrule
       PCA &9.26$\cdot10^{-3}$ & 2.98$\cdot10^{-3}$ \\
         \bottomrule
    \end{tabular}
    \caption{Mean relative test errors, averaged over the three training--test splits}
    \label{tab:errors}
\end{table}

\subsection{Runtime Discussion}
 
In this subsection, we compare the runtimes of the different methods in \Cref{tab:runtimes}. 
We compare the time needed for the computation of the feature map added to the
feature extraction time on the test set (first column),  the time needed for LOOCV, the training time of the final model and the evaluation time on the test set. Note that the full oder model solves and the morphological feature extraction are performed on a different machine (workstation with two Intel Xeon Gold 6240R (48 cores in total)) than the PCA feature extraction and kernel model training (performed on a computer with 64 GB RAM and a 13th Gen Intel i7-13700K processor) due to licensing limitations.
We observe that the LOOCV for the two-layered models is way more expensive than for the shallow models which is due to the optimization of the first layer. 
However, it is still considerably lower (about 30 times for MF-2L-kernel and 7 times for PCA-2L-kernel) than the generation of a single training sample (7514 wcs). 
Furthermore, the LOOCV time and final model training of the MF-kernel models is less expensive than the training on the morphological features. However, 
the computation of the PCA feature map is less expensive than the computation of the morphological features.  
Moreover, the evaluation time of the PCA-kernel models is larger than the evaluation of the MF kernel models. Since the computation of the morphological features i.e.,\ the evaluation of $\Phi_\text{MF}$ takes about 2 wcs per sample, the MF approaches have a larger overall evaluation time (summing evaluation time of MF extraction and MF-1L-kernel/MF-2L-kernel). 
Nevertheless, in both cases the evaluation times are considerably lower than the computation of the training samples. 

\begin{table}[t]
    \centering
    \begin{tabular}{ccccc} 
    \toprule
 model & feat. map & LOOCV time & final model train time & eval. time \\
 \midrule
 Full order model & - & - & - & 9.02$\cdot10^4$ \\
 MF extraction & 9.40$\cdot10^1$ & - & - & 2.40$\cdot10^1$ \\ \hline
MF-1L-kernel & - & 5.90$\cdot 10^0$ & 9.95$\cdot10^{-4}$ & 2.91$\cdot10^{-4}$ \\
MF-2L-kernel & - & 2.01$\cdot10^2$ & 2.12$\cdot10^{-1}$ & 2.67$\cdot10^{-4}$ \\
PCA feat. map comp. & 1.28$\cdot10^1$ & - & - & - \\
PCA-1L-kernel & - & 2.29$\cdot10^1$ & 1.36$\cdot10^{-3}$ & 1.83$\cdot10^{-1}$ \\
PCA-2L-kernel & - & 9.18$\cdot10^3$ & 8.63$\cdot10^{-1}$ & 2.13$\cdot10^{-1}$ \\
\bottomrule
\end{tabular}
    \caption{Runtime comparison of the machine learning methods in wcs.}
    \label{tab:runtimes}
\end{table}

\section{Conclusion and Outlook}\label{sect:concl}
In this work we demonstrated how breakthrough curves can be predicted from the geometry of a 3D porous medium. 
We presented two approaches on how to treat the high dimension of the input data: 
For the first approach, we computed morphological features of the geometries and learned a mapping from these morphological features to the breakthrough curves. 
For the second approach, we computed PCA features of the geometry data set and learned a mapping from these PCA features to the breakthrough curve.

We observed that the MF approach worked well in combination with one-layered kernels and that both approaches worked well in combination with two-layered kernels. 
This is compelling, as the morphological features used in the study are well-known informative and predictive quantities for porous media and in contrast the PCA feature extraction approach is purely data-based.  
Thus, the results indicate the strength of using two-layered kernels, as they automatically adapt to the relevant features. We further observed that the scarce-data situation did not prevent these approaches from predicting well the high-dimensional outputs.

Future work will focus on how ideas from convolutional neural networks can be combined with the framework of deep multi-layered kernels. 
Moreover, we will investigate, whether data-augmentation (rotating/reflecting geometry samples without changing the breakthrough curves) can further improve the PCA-feature approach while avoiding to compute expensive FOM solutions.

\section*{Data availability} 
\begin{minipage}{\linewidth}
The original data presented in the study are openly available at \href{https://doi.org/10.18419/darus-4227}{https://doi.org/10.18419/darus-4227}.
\end{minipage}

\section*{Acknowledgements} 
Funded by BMBF under the contracts 05M20VSA and 05M20AMD (ML-MORE).
The authors acknowledge the funding of the project by the Deutsche Forschungsgemeinschaft (DFG, German Research Foundation) under
Germany's Excellence Strategy - EXC 2075 - 390740016.

\bibliographystyle{plainurl}

\newpage
\begin{appendices}
\renewcommand{\thesection}{\Roman{section}}
\section[]{Geometry Parameters}\label{sect:appendix_S}
In the following, morphological features, i.e.\ geometry parameters are listed,
which were computed based on voxel representation of the porous samples.
We assume that we already have a voxel representation of the domain $\Omega$ and its subsets $\Omega_f,\Omega_w,\Omega_s$, where $\Omega_f,\Omega_w,\Omega_s$ are the free pores, washcoat and solid domains respectively.
\begin{enumerate}
\item Porosity $\epsilon$:
\begin{equation*}
\epsilon \coloneq \frac{|\Omega_{f}|}{|\Omega|}.
\end{equation*}
\item Washcoat volume fraction $\epsilon_{w}$:
\begin{equation*}
\epsilon_w \coloneq \frac{|\Omega_{p}|}{|\Omega|}.
\end{equation*}
\item Volume of the free pores $V$:
\begin{equation*}
V \coloneq |\Omega_f|.
\end{equation*}
\item Surface area for the free pores phase $S$: 
\begin{equation*}
S \coloneq |\partial \Omega_f|.
\end{equation*}
\item Integral of mean curvature for free pores phase, $c_f$. 
For smooth surfaces this quantity is usually defined by the integral
\begin{equation*}
c_f \coloneq \frac{1}{2} \int_{\partial \Omega_f} \left(\frac{1}{k_1}+\frac{1}{k_2} \right) ds,
\end{equation*}
here $ds$ is a surface element of $\partial \Omega_f$ and $k_1$ and $k_2$ are the two principal curvatures from the respective surface element. 
Since the boundaries of the voxelized geometry are piece-wise flat, the software package computes an approximation of this quantity for our phase boundaries.

\item Integral of total curvature for free pores phase, $ct_{f}$. 
For smooth surfaces this quantity is usually defined by the integral
\begin{equation*}
ct_f \coloneq \int_{\partial \Omega_f} \left(\frac{1}{k_1 k_2} \right) ds,
\end{equation*}
here $ds$ is a surface element of $\partial \Omega_f$ and $k_1$ and $k_2$ are the two principal curvatures from the respective surface element. 
The software package computes an approximation of this quantity for our piece-wise flat phase boundaries.
\end{enumerate}              
\end{appendices}

\end{document}